\documentclass[reqno,10pt]{amsart}

\usepackage{xypic}
\input xy
\xyoption{all}
\usepackage{epsfig}
\usepackage{color}
\usepackage{amsthm}
\usepackage{amssymb}
\usepackage{amsmath}
\usepackage{amscd}
\usepackage{amsopn}
\usepackage{url}

\usepackage{hyperref}\hypersetup{colorlinks}

\usepackage{color} 

\definecolor{darkred}{rgb}{1,0,0} 
\definecolor{darkgreen}{rgb}{0,0.8,0}
\definecolor{darkblue}{rgb}{0,0,1}

\hypersetup{colorlinks,
linkcolor=darkblue,
filecolor=darkgreen,
urlcolor=darkred,
citecolor=darkgreen}

\newcommand{\labell}[1] {\label{#1}}

\numberwithin{equation}{section}
\newtheorem {Theorem}{Theorem}
\numberwithin{Theorem}{section}

\newtheorem {Lemma}[Theorem]{Lemma}

\newtheorem {Proposition}[Theorem]{Proposition}

\theoremstyle{definition}
\newtheorem{Definition}[Theorem]{Definition}
\theoremstyle{remark}
\newtheorem{Remark}[Theorem]{Remark}
\newtheorem{Example}[Theorem]{Example}

\expandafter\chardef\csname pre amssym.def at\endcsname=\the\catcode`\@
\catcode`\@=11
\def\undefine#1{\let#1\undefined}
\def\newsymbol#1#2#3#4#5{\let\next@\relax
 \ifnum#2=\@ne\let\next@\msafam@\else
 \ifnum#2=\tw@\let\next@\msbfam@\fi\fi
 \mathchardef#1="#3\next@#4#5}
\def\mathhexbox@#1#2#3{\relax
 \ifmmode\mathpalette{}{\m@th\mathchar"#1#2#3}%
 \else\leavevmode\hbox{$\m@th\mathchar"#1#2#3$}\fi}
\def\hexnumber@#1{\ifcase#1 0\or 1\or 2\or 3\or 4\or 5\or 6\or 7\or 8\or
 9\or A\or B\or C\or D\or E\or F\fi}

\font\teneufm=eufm10
\font\seveneufm=eufm7
\font\fiveeufm=eufm5
\newfam\eufmfam
\textfont\eufmfam=\teneufm
\scriptfont\eufmfam=\seveneufm
\scriptscriptfont\eufmfam=\fiveeufm

\catcode`\@=\csname pre amssym.def at\endcsname

\newcommand{\AC}{{\mathcal A}}

\newcommand{\FF}{{\mathcal F}}

\newcommand{\HH}{{\mathcal H}}

\newcommand{\SC}{{\mathcal S}}

\newcommand{\tH}{\tilde{H}}
\newcommand{\tG}{\tilde{G}}
\newcommand{\tK}{\tilde{K}}
\newcommand{\tx}{\tilde{x}}
\newcommand{\ty}{\tilde{y}}

\def    \R      {{\mathbb R}}
\def    \Z      {{\mathbb Z}}

\newcommand{\const}{{\mathit const}}

\def    \eps      {\epsilon}

\def    \p        {\partial}
\def    \Area   {\operatorname{Area}}

\def    \sgn       {\operatorname{sgn}\,}

\def    \Sp     {\operatorname{Sp}}

\def \codim {\operatorname{codim}} 

\def    \H   {\operatorname{H}} 
\def    \HF   {\operatorname{HF}} 
\def    \CF   {\operatorname{CF}} 
\def    \CZ   {\operatorname{\mu_{\scriptscriptstyle{CZ}}}} 
\def    \CZ  {\operatorname{\mu_{\scriptscriptstyle{CZ}}}} 
\def    \cf   {\operatorname{c}}
\def    \Ho    {\operatorname{\scriptscriptstyle{H}}} 

\def \wM {\omega_{\scriptscriptstyle{M}}}

\begin{document}


\setlength{\smallskipamount}{6pt}
\setlength{\medskipamount}{10pt}
\setlength{\bigskipamount}{16pt}




\title[Coisotropic Maslov Class Rigidity]{On the Maslov Class Rigidity
  for Coisotropic Submanifolds}

\author[Viktor Ginzburg]{Viktor L. Ginzburg}

\address{Department of Mathematics, University of California, Santa Cruz,
  CA 95064, USA} \email{ginzburg@math.ucsc.edu}

\date{\today} 

\thanks{The work is partially supported by the NSF and by the faculty
research funds of the University of California, Santa Cruz.}


\begin{abstract}
  We define the Maslov index of a loop tangent to the characteristic
  foliation of a coisotropic submanifold as the mean Conley--Zehnder
  index of a path in the group of linear symplectic transformations,
  incorporating the ``rotation'' of the tangent space of the leaf --
  this is the standard Lagrangian counterpart -- and the holonomy of
  the characteristic foliation. Furthermore, we show that, with this
  definition, the Maslov class rigidity extends to the class of the
  so-called stable coisotropic submanifolds including Lagrangian tori
  and stable hypersurfaces.
\end{abstract}

\maketitle

\tableofcontents

\section{Introduction and main results }
\labell{sec:main-results}

\subsection{Introduction}
\label{sec:intro}
As the title indicates, the main theme of the paper is the Maslov class
rigidity for coisotropic submanifolds. To be more specific, we define
the Maslov index of a loop tangent to the characteristic foliation in
a coisotropic submanifold and show that a displaceable, stable
coisotropic submanifold carries a loop with Maslov index in the range
$[1,\, 2n+1-k]$, where $2n$ is the dimension of the ambient manifold
and $k$ is the codimension of the coisotropic submanifold.

The study of symplectic topology of coisotropic submanifolds can be
traced back to Moser's paper \cite{Mo} followed by \cite{Ba,EH,Ho90a}
and by the work of Bolle, \cite{Bo96,Bo98}. Recently, the field has
entered a particularly active phase; see
\cite{AF08,AF08g,Dr,Gi07coiso,Gu09imrn,Gu:ex,Ka,Ke07coiso,To,U09coiso,Zilt08,ZiltPrep}.
Most of these papers, with the exception of \cite{ZiltPrep}, concern
such questions as generalizations to coisotropic submanifolds of the
Lagrangian intersection property or of the existence of closed
characteristics on stable hypersurfaces. The present work, which can
be thought of as a follow-up to \cite{Gi07coiso}, focuses mainly on
the coisotropic version of the Maslov class rigidity, also considered
in \cite{ZiltPrep}.

The aspect of the Maslov class rigidity we are concerned with here is
the fact that the Maslov class of a closed displaceable Lagrangian
submanifold automatically satisfies certain restrictions. Namely, the
minimal Maslov number of such a submanifold lies between $1$ and
$n+1$. This phenomenon was originally studied in
\cite{Po91Tr,Po91MZ,Vi90} and there are two methods of proving 
results of this type. One of these methods uses the holomorphic curves
technique (see \cite{AL,Po91Tr,Po91MZ}) and at this moment it is not
known how to directly apply it to coisotropic submanifolds due to the
lack of Fredholm properties for the Cauchy--Riemann problem with
coisotropic boundary conditions. The second approach, originating from
\cite{Vi90}, relies on Hamiltonian Floer homology (or its equivalent)
and in combination with certain estimates from \cite{Bo98} can be
easily adapted to the coisotropic setting; see, e.g.,
\cite{Gi07coiso}. Here, we heavily draw from the modern interpretation
of this method given in \cite{Ke09,KS}.

The Maslov index of a loop tangent to the characteristic foliation is
the mean Conley--Zehnder index of a certain path in $\Sp(2n)$
associated with the loop and comprising the ``rotation'' of the
tangent space of the leaf, as the standard Lagrangian counterpart, and
the holonomy of the characteristic foliation. Hence, the index can be
an arbitrary real number. This definition, which can also be found in
\cite{ZiltPrep} where it is treated in great detail, is of independent
interest.  Then, the proof of the Maslov class rigidity for
coisotropic submanifolds follows the path of \cite{Ke09,KS,Vi90}.  The
main new element of the proof is that we circumvent relating the
Conley--Zehnder and Morse indices as in \cite{Du,Vi90}, but instead
use the explicit expression for the geodesic flow of a metric,
capitalizing on the fact that the submanifolds in question are stable
and hence admit a leaf-wise flat metric.

\subsection{Coisotropic Maslov index}
\label{sec:maslov-def}
Let $M$ be a coisotropic submanifold of a symplectic manifold
$(W^{2n},\omega)$. Denote by $\FF$ the characteristic foliation of
$M$; see Section \ref{subsec:prelim} for the definition. The normal
bundle $T^\perp M$ to $M$ is canonically isomorphic to the (leaf-wise)
cotangent bundle $T^*\FF$ to $\FF$ and the direct sum $T\FF\oplus
T^\perp M$ is a symplectic vector bundle over $M$. Furthermore, we
have a symplectic vector bundle decomposition
\begin{equation}
\label{eq:split}
TW\mid_M=(T\FF\oplus T^\perp M)\oplus T^\perp\FF,
\end{equation}
where $T^\perp\FF$ is the normal bundle to $\FF$ in $M$. Note that
$T^\perp\FF$ carries a symplectic leaf-wise flat connection.

Consider a loop $\gamma\colon S^1\to M$ tangent to $\FF$, contractible
in $W$ and equipped with a capping $u\colon D^2\to W$. The capping $u$
gives rise to a symplectic trivialization $\zeta$, unique up to
homotopy, of the pull-back bundle $\gamma^* TW$.  Let us assume first
that $T\FF$ is orientable along $\gamma$ (i.e., the pull-back
$\gamma^* T\FF$ is orientable), and hence trivial, and fix a
trivialization $\xi$ of this vector bundle. Then the pull-back
$\gamma^*(T\FF\oplus T^\perp M)$ receives a symplectic trivialization
$\xi\oplus \xi^*$. This trivialization can be viewed as a family of
symplectic maps $\Xi(t)\colon T_{\gamma(0)}\FF\oplus
T^\perp_{\gamma(0)}M\to T_{\gamma(t)}\FF\oplus T^\perp_{\gamma(t)}M$
parametrized by $t\in S^1$.  Combining the family $\Xi(t)$ with the
holonomy $\Gamma(t)\colon T^\perp_{\gamma(0)} \FF\to
T^\perp_{\gamma(t)} \FF $ along $\gamma$, we obtain a family of
symplectic maps $\Xi(t)\oplus \Gamma(t)\colon T_{\gamma(0)}W\to
T_{\gamma(t)}W$, which, using the trivialization $\zeta$, we can
regard as a path $\Phi\colon [0,\,1]\to\Sp(2n)$.

\begin{Definition}
\label{def:maslov}
The \emph{coisotropic Maslov index} $\mu(\gamma,u)$ of the capped loop
$(\gamma, u)$ is the negative mean Conley--Zehnder index $-\Delta(\Phi)\in\R$. (We
refer the reader to \cite{Lo,SZ} for a detailed discussion of the mean
index; here we use the notation and conventions from \cite{GG09wm};
see Section \ref{subsec:conv}.) 
When $T\FF$ is not orientable along
$\gamma$, we set $\mu(\gamma,u):=\mu(\gamma^2,u^2)/2$, where
$(\gamma^2,u^2)$ stands for the double cover of $(\gamma,u)$.
\end{Definition}

The standard argument shows that the index $\mu(\gamma,u)$ is well
defined, i.e., independent of the choice of the trivializations $\xi$
and $\zeta$. Furthermore, it is also independent of the choice of
splitting \eqref{eq:split}: the normal bundle $T^\perp\FF$ is
unambiguously defined only as the quotient $TW/T\FF$ while the
splitting requires a choice of the complement to $T\FF$ in $TW$. To
see that $\Delta(\Phi)$ is independent of this choice, we argue as
follows; cf.\ the proof of \cite[Lemma 2.6]{GG09wm}. Observe that the
path $\tilde{\Phi}$ resulting from a different splitting is homotopic
to the concatenation of the path $\Phi$ with a path $\Psi$ of the form
$\Psi(t)=I+A(t)$, where $I$ is the identity map and $A(t)\colon
T^\perp\FF\to (T\FF\oplus T^\perp M)$. Thus, all eigenvalues of
$\Psi(t)$ are equal to one and, as a consequence,
$\Delta(\Psi)=0$. Hence, by the additivity and homotopy invariance of
the mean index (see, e.g., \cite{GG09wm,Lo,SZ}), we have
$\Delta(\tilde{\Phi})=\Delta(\Phi)$. 

It is worth emphasizing that, in contrast with the ordinary Lagrangian
Maslov index, the coisotropic Maslov index is not, in general, an
integer and that this index is different from the one considered in
\cite{Oh}.  The negative sign in the definition of the coisotropic
Maslov index is, of course, a matter of conventions: this is the price
we have to pay to match the sign of the standard Maslov index for
Lagrangian submanifolds (Example \ref{ex:Lagr}) while using the
conventions from \cite{GG09wm}; see Section~\ref{subsec:conv}.

It is easy to see that the coisotropic Maslov index has the following
properties.
\begin{itemize}
\item Homotopy invariance: $\mu(\gamma,u)$ is invariant, in the obvious sense, under a
  homotopy of $\gamma$ in a leaf of $\FF$. In particular,
  $\mu(\gamma,u)=0$ when $u$ is homotopic (rel boundary) to a disc in
  the leaf of $\FF$ containing $\gamma$.

\item Recapping: $\mu(\gamma,u\# v)=\mu(\gamma,u)-2\left<c_1(TW),v\right>$, where
  the capping $u\# v$ is obtained by attaching the sphere
  $v\in \pi_2(W)$ to $u$. In particular, $\mu(\gamma):=\mu(\gamma,u)$ is
  independent of $u$ when $c_1(TW)\mid_{\pi_2(W)}=0$.

\item Homogeneity: $\mu(\gamma^k,u^k)=k\mu(\gamma,u)$, where
  $(\gamma^k,u^k)$ stands for the $k$-fold cover of
  $(\gamma,u)$. Moreover, when $c_1(TW)\mid_{\pi_2(W)}=0$, the Maslov
  index gives rise to a homomorphism $\pi_1(F)\to\R$ for any leaf $F$
  of $\FF$.

\end{itemize}

\begin{Example} 
\label{ex:Lagr}
When $M$ is a Lagrangian submanifold of $W$, the foliation $\FF$ has
only one leaf, the manifold $M$ itself, and the coisotropic Maslov
index coincides with the ordinary Maslov index. Indeed, in this case, Definition
\ref{def:maslov} turns into one of the definitions of this index.
\end{Example}

\begin{Example} 
\label{ex:hyper}
When $u$ is contained in $M$, the index
  $\mu(\gamma,u)$ is equal to the mean index of the holonomy along
  $\gamma$ with respect to a symplectic trivialization of $T^\perp
  \FF$ associated with $u$. For instance, when $M$ is a regular level
  of a Hamiltonian and $\gamma$ is a periodic orbit (and again $u$
  is contained in $M$), the Maslov index $\mu(\gamma,u)$ is equal to
  the mean index of $\gamma$ in $M$.
\end{Example}

\begin{Example} When all leaves of $\FF$ are closed and form a
  fibration, the path $\Phi$ is a loop and $\mu(\gamma,u)$ is equal to
  the Maslov index of this loop.  (In particular, then $\mu(\gamma,u)$
  is an integer.) In this setting, the coisotropic Maslov index is further
  investigated by Ziltener, \cite{ZiltPrep}. Moreover, one can
  express the coisotropic Maslov index via the Lagrangian Maslov index
  in the graph of $\FF$; see \cite{Zilt08,ZiltPrep} for details.
\end{Example}

Now we are in a position to state the main result of the paper. A much
more detailed discussion of the coisotropic Maslov index can be found
in \cite{ZiltPrep}.

\subsection{Rigidity of the coisotropic Maslov index}
\label{sec:rigidity}
Let $W$ be a symplectically aspherical manifold, which we assume to be
either closed or geometrically bounded and wide (e.g., convex at
infinity) in the sense of \cite{Gu07}.

\begin{Theorem}
\label{thm:main}
Let $W^{2n}$ be as above and let $M^{2n-k}\subset W$ be a closed,
stable, displaceable coisotropic submanifold. (See Section
\ref{subsec:prelim} for the definitions.) Then, for any $\delta>0$,
there exists a loop $\eta$ tangent to $\FF$ and contractible in $W$
and such that
\begin{equation}
\label{eq:bound-index}
1\leq \mu(\eta)\leq 2n+1-k
\end{equation}
and 
\begin{equation}
\label{eq:bound-area}
0<\Area(\eta)\leq e(M)+\delta,
\end{equation}
where $\Area(\eta)$ is the symplectic area bounded by $\eta$ and
$e(M)$ is the displacement energy of $M$.
\end{Theorem}

\begin{Example}
  As in Example \ref{ex:Lagr}, assume that $M$ is a stable Lagrangian
  submanifold (and hence a torus). Then $k=n$ and the theorem reduces
  to a particular case of the standard Lagrangian Maslov class
  rigidity.  This version of rigidity is established in \cite{Vi90} for
  $W=\R^{2n}$ and in \cite{Ke09,KS} for closed ambient manifolds; see
  also \cite{AL,Po91Tr,Po91MZ} for generalizations.
\end{Example}

\begin{Example}
  Assume that $M$ is a stable, displaceable, simply connected
  hypersurface. Then, by \eqref{eq:bound-index} and Example
  \ref{ex:hyper}, $M$ carries a closed characteristic $\eta$ with
  $1\leq \Delta(\gamma)\leq 2n$.  This is apparently a new
  observation. However, if we replace the upper bound by $2n+1$, the
  assertion becomes an easy consequence of the properties of the
  mean index and, for instance, the displacement or symplectic homology
  proof of the almost existence theorem; see, e.g., \cite{FHW,Gi:alan,Gu07,HZ94} and
  references therein.
\end{Example}

\begin{Remark} A word on the hypotheses of the theorem is due
  now. The assumption that $W$ be symplectically aspherical is imposed
  here only for the sake of simplicity and can be significantly
  relaxed along the lines of
  \cite{Ke07coiso,U09coiso}. Hypothetically, a combination of our
  argument with the reasoning from these works should lead to a
  generalization of the theorem to the case where we only require the
  subgroup $\left<\omega,\pi_2(M)\right>\subset \R$ to be discrete as
  in \cite[Theorem 1.6]{U09coiso} or, at least, where $W$ is monotone
  or negative monotone; see \cite{Ke07coiso}. (In such a
  generalization, the geodesic $\eta$ is, of course, equipped with
  capping.)

  The condition that $M$ is stable cannot be entirely omitted due to the
  counterexamples to the Hamiltonian Seifert conjecture showing that
  there exist hypersurfaces in $\R^{2n}$ ($C^2$ when $2n=4$) without
  closed characteristics; see \cite{Gi99:orbits,GG03:seifert} and
  references therein. However, this condition can possibly be relaxed
  as in \cite[Section~7]{U09coiso}.

  Finally note that the existence of a loop $\eta$ satisfying
  \eqref{eq:bound-area} is established in
  \cite[Theorem~2.7]{Gi07coiso}, where the second inequality (with
  $\delta=0$) is proved under the additional hypothesis that $M$ has
  restricted contact type. Thus, even when only the area bounds are
  concerned, Theorem \ref{thm:main} is a generalization (up to the
  issue of $\delta$) of the results from \cite{Gi07coiso}, which
  became possible due to incorporating a technique from \cite{Ke09,KS}
  into the proof.
\end{Remark}

\begin{Remark} 
  It is tempting to conjecture that the Maslov class of $M$ is still
  non-zero even when the stability assumption in Theorem
  \ref{thm:main} is dropped and all leaves of $\FF$ may be
  contractible. However, it is not entirely clear how to define
  this Maslov class and what cohomology space this class should lie
  in. The situation contrasts sharply with a similar question for the
  Liouville class of $M$, which can always be defined, when $W$ is
  exact, as the class $[\lambda|_\FF]$ of a global primitive $\lambda$
  of $\omega$ in the tangential de Rham cohomology $\H^1(\FF)$; see
  \cite[Section 1.2]{Gi07coiso}.
\end{Remark}

\subsection*{Acknowledgments} 
The author is grateful to Ba\c sak G\"urel and Ely Kerman for 
useful discussions and remarks. He would also like to
thank Yael Karshon and Fabian Ziltener for an inspiring discussion
of the notion of the coisotropic Maslov index.


\section{Preliminaries}
\label{sec:lw-int}
We start this section by recalling the relevant definitions and basic
results concerning coisotropic submanifolds.  In Section
\ref{subsec:conv}, we set our conventions and notation.  

\subsection{Stable coisotropic submanifolds}
\label{subsec:prelim}
Let, as above, $(W^{2n}, \omega)$ be a symplectic manifold and let
$M\subset W$ be a closed, coisotropic submanifold of codimension
$k$. Set $\wM=\omega|_M$. Then, as is well known, the distribution
$\ker \wM$ has dimension $k$ and is integrable. Denote by $\FF$ the
characteristic foliation on $M$, i.e., the $k$-dimensional foliation
whose leaves are tangent to the distribution $\ker \wM$.

\begin{Definition}
\label{def:B}
The coisotropic submanifold $M$ is said to be \emph{stable} if there
exist one-forms $\alpha_1,\ldots,\alpha_k$ on $M$ such that $\ker
d\alpha_i\supset \ker\wM$ for all $i=1,\ldots,k$ and
\begin{equation}
\labell{eq:ct}
\alpha_1\wedge\cdots\wedge\alpha_k\wedge \wM^{n-k}\neq 0
\end{equation}
anywhere on $M$. We say that $M$ has \emph{contact type} if the forms
$\alpha_i$ can be taken to be primitives of $\wM$. Furthermore, $M$
has \emph{restricted} contact type if the forms $\alpha_i$ extend to
global primitives of $\omega$ on $W$.
\end{Definition}

Stable and contact type coisotropic submanifolds were introduced by
Bolle in \cite{Bo96,Bo98} and considered in a more general setting in
\cite{Gi07coiso} and also by Kerman, \cite{Ke07coiso}, and Usher,
\cite{U09coiso}.  We refer the reader to \cite{Gi07coiso} for a
discussion of the requirements of Definition \ref{def:B} and
examples. Here we only note that although Definition \ref{def:B} is
natural, it is quite restrictive. For example, a stable Lagrangian
submanifold is necessarily a torus and a stable coisotropic
submanifold is automatically orientable.

Assume henceforth that $M$ is stable. Then the normal bundle $T^\perp
M$ to $M$ in $W$ is trivial, since it is isomorphic to $T^*\FF$ and
the latter bundle is trivial due to \eqref{eq:ct}.  From now on, we
fix the trivialization $T^\perp M=T^*\FF\cong M \times
\R^k$ given by the forms $\alpha_i$ and identify a small neighborhood
of $M$ in $W$ with a neighborhood of $M$ in $T^*\FF = M \times
\R^k$. We will use the same symbols $\wM$ and $\alpha_i$ for
differential forms on $M$ and for their pullbacks to $M \times
\R^k$. (In other words, we are suppressing the pullback notation
$\pi^*$, where $\pi \colon M \times \R^k \to M$, unless its presence
is absolutely necessary.)  As a consequence of the Weinstein
symplectic neighborhood theorem, we have

\begin{Proposition}[\cite{Bo96, Bo98}]
\label{prop:normalform}
Let $M$ be a closed, stable coisotropic submanifold of $(W^{2n},
\omega)$ with $ \codim M = k$.  Then, for a sufficiently small $r>0$,
there exists a neighborhood of $M$ in $W$, which is symplectomorphic to
$U_r=\{ (q,p) \in M \times \R^k \mid |p| < r \}$ equipped with the
symplectic form $\omega = \wM + \sum_{j=1}^k d(p_j \alpha_j)$. Here
$(p_1,\ldots,p_k)$ are the coordinates on $\R^k$ and $|p|$ is the
Euclidean norm of $p$.
\end{Proposition}

Thus, a neighborhood of $M$ in $W$ is foliated by a family of
coisotropic submanifolds $M_p = M\times \{p\} $ with $p \in B^k_r$,
where $B_r^k$ is the ball of radius $r$ centered at the origin in
$\R^k$. Moreover, a leaf of the characteristic foliation on $M_p$
projects onto a leaf of the characteristic foliation on $M$.

Furthermore, we have

\begin{Proposition}[\cite{Bo96, Bo98, Gi07coiso}]
\label{prop:flow}
Let $M$ be a stable coisotropic submanifold.

\begin{enumerate}

\item[(i)] The leaf-wise metric $(\alpha_1)^2+\cdots+(\alpha_k)^2$ on
  $\FF$ is leaf-wise flat.

\item[(ii)] The Hamiltonian flow of
  $\rho=(p_1^2+\cdots+p_k^2)/2=|p|^2/2$ is the leaf-wise geodesic flow
  of this metric.

\end{enumerate}

\end{Proposition}

We conclude this section by pointing out that the metric $\rho$
extends to a true metric on $M$ such that the leaves of $\FF$ are
totally geodesic submanifolds and that the existence of such a metric
is equivalent to the stability of $M$ when $M$ is a hypersurface; see
\cite{Su} and \cite[Section 7]{U09coiso}.

\subsection{Conventions and notation} 
\label{subsec:conv}
In this section we specify conventions and notation used throughout the paper. 

\subsubsection{Action functional and the Hamilton equation}
Let $(W^{2n},\omega)$ be a symplectically aspherical
manifold, i.e., $\omega|_{\pi_2(W)}=c_1|_{\pi_2(W)}=0$.  Denote by
$\Lambda W$ the space of smooth contractible loops $\gamma\colon
S^1\to W$ and consider a time-dependent Hamiltonian $H\colon S^1\times
W\to \R$, where $S^1=\R/\Z$.  Setting $H_t = H(t,\cdot)$ for $t\in
S^1$, we define the action functional $\AC_H\colon \Lambda W\to \R$ by
$$
\AC_H(\gamma)=\AC(\gamma)+\int_{S^1} H_t(\gamma(t))\,dt,
$$
where $\AC(\gamma)=-\Area(\gamma)$ is the negative symplectic area
bounded by $\gamma$. In other words,
$$
\AC(\gamma)=-\int_{u}\omega,
$$
where $u\colon D^2\to W$ is a capping of $\gamma$, i.e., 
$u|_{S^1}=\gamma$.  The least action principle asserts that
the critical points of $\AC_H$ are exactly the contractible one-periodic
orbits of the time-dependent Hamiltonian flow $\varphi_H^t$ of $H$,
where the Hamiltonian vector field $X_H$ of $H$ is defined by the Hamilton equation
$i_{X_H}\omega=-dH$. 

\subsubsection{Conley--Zehnder  index} 
Consider a finite--dimensional symplectic vector space $V$.  We denote
by $\Sp(V)$ the group of linear symplectic transformations of $V$ and,
as usual, set $\Sp(2n)=\Sp(\R^{2n})$. Furthermore, we let
$\Delta(\Phi)$ stand for the mean index of a path $\Phi\colon
[0,\,T]\to\Sp(V)$ and, when $\Phi$ is non-degenerate (i.e., $\Phi(T)$
has no eigenvalues equal to one), we denote by $\CZ(\Phi)$ the
Conley--Zehnder index of $\Phi$. We refer the reader to
\cite{Lo,Sa,SZ} and also \cite{GG09wm} for the definitions and a
detailed discussion of these notions. In this paper, we normalize
these indices as in \cite{GG09wm}. This normalization is different
from the ones in \cite{Lo,Sa,SZ}. For instance, our $\CZ(\Phi)$ is the
negative of the Conley--Zehnder index as defined in \cite{Sa}. For the
flow $\Phi(t)$ with $0\leq t\leq 1$ generated by a non-degenerate
quadratic Hamiltonian $H$ with small eigenvalues, we have
$\CZ(\Phi)=-\sgn (H)/2$, where $\sgn(H)$ is the signature of $H$ (the
number of positive squares minus the number of negative squares).  In
particular, when $H$ is negative definite, we have $\CZ(\Phi)=n$ where
$2n=\dim V$ and $\Delta(\Phi)>0$. In other words, when $\CZ(\Phi)$ is
interpreted as the intersection index of $\Phi$ with the discriminant
$\Sigma\subset \Sp(V)$ formed by symplectic transformations with at
least one eigenvalue equal to one, $\Sigma$ is co-oriented by the
Hamiltonian vector field of a negative definite Hamiltonian. 

Recall also from \cite{SZ} that, regardless of conventions, we have
\begin{equation}
\label{eq:mean-index}
|\Delta(\Phi)-\CZ(\Phi)|< n \text{ and }
\Delta(\Phi)=\lim_{k\to\infty}\frac{\CZ(\Phi^k)}{k},
\end{equation}
where in the inequality we require
$\Phi(T)$ to be non-degenerate and, in the limit identity, we 
assume that $\Phi(T)^k\not \in \Sigma$ for all $k$ and thus
$\CZ(\Phi^k)$ is defined. Note that here we can replace $\Phi^k$ by the
concatenation of the paths $\Phi$, $\Phi(T)\Phi$, etc, up to
$\Phi(T)^{k-1}\Phi$.

Let now $x$ be a contractible periodic orbit of $H$ on $W^{2n}$. Using a
trivialization of $x^*TW$ arising from a capping of $x$, we can
interpret the linearized flow $d\varphi_H^t$ along $x$ as a path
$\Phi$ in $\Sp(2n)$. The mean index $\Delta(x)$ of $x$ is by definition
$\Delta(\Phi)$. When $x$ is non-degenerate, we also set
$\CZ(x):=\CZ(\Phi)$. Since $c_1(TW)|_{\pi_2(W)}=0$, these
indices are well-defined, i.e., independent of the capping.  When we need
to emphasize the role of $H$, we write $\Delta_H(x)$ and
$\CZ(x,H)$. By \eqref{eq:mean-index}, we have
\begin{equation}
\label{eq:mean-index2}
|\Delta(x)-\CZ(x)|< n \text{ and }
\Delta(x)=\lim_{k\to\infty}\frac{\CZ(x^k)}{k}.
\end{equation}
As in \eqref{eq:mean-index}, we require here  $x$ to be
non-degenerate for $\CZ(x)$ to be defined,  and,  in the limit identity, we
assume that $x$ is strongly non-degenerate, i.e., all iterated
orbits $x^k$ are non-degenerate.  Finally note that with our
normalizations $\Delta(x)>0$ and $\CZ(x)=n$ when $x$ is a
non-degenerate maximum (with small Hessian) of an autonomous
Hamiltonian.

\subsubsection{Floer homology}
In the definition of Floer homology, we adopt literally the
conventions and notation from \cite{Gi07coiso}. All Hamiltonians
considered in this paper are assumed to be compactly supported. The
manifold $W$, in addition to being symplectically aspherical, is
required to be either closed or geometrically bounded and wide in the
sense of \cite{Gu07}.  (See, e.g., \cite{AL,CGK,Si} for the precise
definition and a discussion of geometrically bounded manifolds.)

Examples of geometrically bounded manifolds include symplectic
manifolds which are convex at infinity (e.g., $\R^{2n}$ and cotangent
bundles) as well as twisted cotangent bundles.  
Under the hypotheses that $W$ is symplectically aspherical and
geometrically bounded, the compactness theorem for Floer's connecting
trajectories holds (see \cite{Si}) and the filtered $\Z$-graded Floer
homology of a compactly supported Hamiltonian on $W$ is defined for
action intervals not containing zero; see, e.g., \cite{CGK,GG04} and
references therein. We use the wideness hypothesis in Section
\ref{sec:pinned} when considering a version of the ``pinned'' action
selector introduced in \cite{Ke09}. This requirement is not
restrictive, for, to the best of the author's knowledge, no examples
of geometrically bounded open manifolds that are not wide are known.

We use the notation $\HF_*^{(a,\,b)}(H)$ for the filtered Floer
homology of $H$, graded by the Conley--Zehnder index.  The end-points
$a$ and $b$ are always assumed to be outside the action spectrum
$\SC(H)$ of $H$ and, if $W$ is open, we require that
$0\not\in(a,\,b)$. When $W$ is closed, we have a canonical isomorphism
$\HF_*(H)=\H_{*+n}(W;\Z_2)$, where as usual
$\HF_*(H)=\HF_*^{(-\infty,\,\infty)}(H)$. When all periodic orbits of
$H$ with action in $(a,\,b)$ are non-degenerate, we let
$\CF_*^{(a,\,b)}(H)$ be the vector space generated over $\Z_2$ by such
orbits, graded by the Conley--Zehnder index. The downward Floer
differential $\p\colon \CF_*^{(a,\,b)}(H)\to \CF_{*-1}^{(a,\,b)}(H)$
is then defined in the standard way and $\HF_*^{(a,\,b)}(H)$ is the
homology of the resulting Floer complex. The above non-degeneracy
requirement is generic (as long as $0\not\in(a,\,b)$ if $W$ is open)
and, in general, we set $\HF_*^{(a,\,b)}(H):=\HF_*^{(a,\,b)}(\tH)$,
where $\tH$ is a small perturbation of $H$ having only non-degenerate
orbits with action in $(a,\,b)$. Since $a$ and $b$ are outside
$\SC(H)$, the homology $\HF_*^{(a,\,b)}(\tH)$ is independent of $\tH$
as long as $\tH$ is sufficiently close to $H$. We refer the
reader to \cite{CGK,Gi07coiso,GG04} for the proofs and further details
on the construction and properties of the Floer homology in this
setting as well as for further references.


\section {Proof of the main theorem}
\label{sec:proof}

\subsection{Maslov index for stable coisotropic submanifolds}
\label{subsec:index}
Let $M$ be a stable coisotropic submanifold.  In this section, we
interpret the mean index $\Delta_\rho(x)$ of a periodic orbit $x$ of
the leaf-wise geodesic flow on $M$ as, up to a sign, the coisotropic
Maslov index of the projection $\gamma$ of $x$ to $M$. Furthermore, we
establish certain bounds, going beyond \eqref{eq:mean-index2}, on the
Conley--Zehnder index of a small non-degenerate perturbation of $x$.
Throughout this subsection, we will use the notation from
Section~\ref{subsec:prelim}. In particular, we fix a neighborhood
$U=M\times B$, where $B=B_r$, of $M$ in $W$.  Thus, let $x$ be a
non-trivial, contractible in $W$ closed orbit of the Hamiltonian flow
of $\rho$ and let $\gamma=\pi(x)$. Then $\gamma$ is also contractible
in $W$.

\begin{Proposition}
\label{prop:index}
We have 
\begin{equation}
\label{eq:index-rho-mu}
\mu(\gamma)=-\Delta_\rho(x). 
\end{equation}
\end{Proposition}

\begin{proof}
  It is convenient to first extend the decomposition \eqref{eq:split}
  from $TW|_M$ to $TW|_U$ as follows. Recall from Section
  \ref{subsec:prelim} that the submanifolds $M_p=M\times\{p\}\subset
  M\times B$, with $p\in B$, are coisotropic and that the
  characteristic foliation $\FF_p$ of $M_p$ projects to $\FF$ under
  $\pi$.  Denote by $\tilde{\FF}$ the resulting foliation of $U$,
  obtained as the union of foliations $\FF_p$.  Furthermore, let $TM$
  be the horizontal tangent bundle in $M\times B$, i.e.,
  $(TM)_{(q,p)}=T_{(q,p)}M_p$ where $(q,p)\in U=M\times B$, and
  likewise let $TB$ denote the vertical bundle $\ker\pi_*$. Then the
  normal bundle $T^\perp \tilde{\FF}$ to $T\tilde{\FF}$ in $TM$ can be
  realized as the sub-bundle $E=(\cap_i\ker\pi^*\alpha_i) \cap TM$. We have
  the symplectic decomposition
\begin{equation}
\label{eq:split2}
TW=(T\tilde{\FF}\oplus TB)\oplus E, 
\end{equation}
which turns into \eqref{eq:split} once restricted to $M$. 

The linearized projection $\pi_*$ gives rise to an isomorphism between
the fibers $(T\tilde{\FF})_{(q,p)}$ and $T_q\FF$, and $E_{(q,p)}$ and
$T^\perp_q\FF$. Furthermore, $(TB)_{(q,p)}$ is naturally isomorphic to
$T_0B=T^\perp_q M$.  Thus, we have a (symplectic) linear isomorphism
between the decomposition \eqref{eq:split2} along $x$ and
\eqref{eq:split} along $\gamma$. In particular, we obtain an
isomorphism between the bundles $x^*TW$ and $\gamma^* TW$ giving rise
to a one-to-one correspondence between trivializations of $TW$ along
$x$ and along $\gamma$. In what follows, we fix a trivialization
arising from a capping of $x$.

Furthermore, recall that the flow of $\rho$ on $U$ can be identified
with the geodesic flow of the leaf-wise metric $\rho$ on $M$.  Thus, we
need to prove that the mean index of the linearized geodesic flow
$G(t)$ along $x$ is equal to $\Delta(\Phi)$.  The geodesic flow
preserves the terms $T\tilde{\FF}\oplus TB$ and $E$ in the
decomposition \eqref{eq:split2}. Indeed, the fact that the first term
is conserved is clear: the geodesic flow is tangent to the leaves. To
show that the second term is conserved, it suffices to recall that, as
mentioned above, the flow is tangent to the manifolds $M_p$ due to
conservation of momenta and that the restrictions
$\pi^*\alpha_j\mid_{M_p}$ are conserved since
$L_{X_\rho}\pi^*\alpha_j=dp_j$.

Next let us show that
\begin{equation}
\label{eq:normal}
G|_E=\Gamma ,
\end{equation}
where we identified $x^*E$ and $\gamma^*T^\perp \FF$.
To this end, let us recall the definition of the holonomy
$\Gamma$. Consider an element $[v]$ in $T_{\gamma(0)}^\perp\FF
=T_{\gamma(0)}M/T_{\gamma(0)}\FF$ represented by a vector $v\in
T_{\gamma(0)}M$. (Here and below, it is more convenient to think of $E$
and $T^\perp\FF$ as quotient bundles rather than
sub-bundles.)  Let $\eta\colon [0,\delta)\to M$ be a smooth map with
$\eta(0)=\gamma(0)$ and $\eta'(0)=v$. Let now $\gamma$ be parametrized
by, say, $[0,\,T]$ and let $\sigma\colon [0,\,T]\times [0,\,\delta)\to
M$ be a map whose restriction to $[0,\,T]\times 0$ is $\gamma$, to
$0\times [0,\,\delta)$ is $\eta$ and such that $\sigma|_{[0,\,T]\times
  s}$, for all $s\in [0,\,\delta)$, lies in a leaf of $\FF$. The class
$[(\p \sigma/ \p s)(t,0)]\in T^\perp_{\gamma(t)}\FF$ is
independent of the choice of $\sigma$ and is the image $\Gamma(t)[v]$.
Let now $w(s)\in T_{\eta(s)}\FF$ be a smooth family of vectors tangent
to $\FF$ and such that $w(0)=\dot{\gamma}(0)$. Consider the
parametrized surface $\sigma$ defined by setting
$\sigma|_{[0,\,T]\times s}$ to be the leaf-wise geodesic with the initial
conditions $(\gamma(s),w(s))$. Then, in particular, 
$[(\p \sigma/ \p s)(t,0)]$ is independent of the choice of the curve
$\eta$ and the family $w$.  Furthermore, on the one hand, this vector
represents $G(t)[v]$ by the definition of the linearized geodesic
flow and, on the other, it is $\Gamma(t)[v]$ due to the above
description of the holonomy.

To complete the argument, it would be sufficient to show that
$G|_{T\tilde{\FF}\oplus TB}=\Xi$, where we identified
$x^*(T\tilde{\FF}\oplus TB)$ and $\gamma^*(T\FF\oplus T^\perp M)$, but
this is not true. Let us fix a basis $\xi(0)\in
T_{\gamma(0)}\FF$. Then, since the metric is flat, $G(t)\xi(0)$ is the
basis $\xi(t)$ in $T_{\gamma(t)}\FF$ obtained form $\xi(0)$ by the
parallel transport along $\gamma$. Let $\xi^*(0)\in
T^*_{\gamma(0)}\FF=T^\perp_{\gamma(0)}\FF$ be the basis dual to
$\xi(0)$. Then $G(t)\xi^*(0)=t\xi(t)+\xi^*(t)\in
T_{\gamma(t)}\FF\oplus T^*_{\gamma(t)}\FF$ in obvious notation. We
conclude that $G(t)|_{T\tilde{\FF}\oplus TB}=\Xi(t) +A(t)$, where
$A(t)\colon T^*_{\gamma(t)}\FF\to T_{\gamma(t)}\FF$.

To finish the proof, we argue as when showing in Section
\ref{sec:maslov-def} that the coisotropic Maslov index is independent
of the splitting \eqref{eq:split}. With a trivialization fixed, we can
view $G$ and $\Phi=\Xi\oplus \Gamma$ as paths in $\Sp(2n)$. Then, $G$
is homotopic with fixed end-points to the concatenation of $\Phi$ and 
the path $\Psi(t)=I+A(t)$. All eigenvalues of
$\Psi(t)$ are equal to one and therefore
$\Delta(\Psi)=0$. Thus, by the additivity and homotopy invariance of
the mean index (see, e.g., \cite{GG09wm,Lo,SZ}), we have
$\Delta(G)=\Delta(\Phi)=:-\mu(\gamma)$.
\end{proof}

\begin{Remark}
  Proposition \ref{prop:index} has the following hypothetical
  generalization. Assume that $M$ admits a metric with respect to
  which $\FF$ is totally geodesic. Referring the reader to
  \cite[Section 7]{U09coiso} for a detailed discussion of this
  condition, we only mention here that it is satisfied when $M$ is
  Lagrangian (for any metric on $M$) and when $M$ is stable. In the
  latter case, $\FF$ is totally geodesic with respect to $\rho$. Then,
  conjecturally, the mean Conley--Zehnder index of $x$ is equal, up
  to a sign, to the sum of the mean Morse index of $\gamma$ and
  $\mu(\gamma)$. When $M$ is stable, the mean Morse index is zero
  since $\rho$ is flat, and this conjecture reduces to Proposition
  \ref{prop:index}. When $M$ is Lagrangian and $x$ is non-degenerate,
  the conjecture essentially reduces to a well-known
  relation between the Conley--Zehnder, Morse, and Maslov indices. The
  latter is proved in \cite{Vi90} using the results from \cite{Du} in
  the context of the finite-dimensional reduction. A proof relying on the
  Floer theory version of the Conley--Zehnder index can be found in, e.g.,
  \cite{We}; see also \cite{KS} for a simple argument.
\end{Remark}

The next proposition is a substitute for the relation between the
Conley--Zehnder and Maslov indices.

\begin{Proposition}
\label{prop:CZ-Delta}
Let $K$ be a small perturbation of $\rho$ and $\tx$ be a non-degenerate
periodic orbit of $K$ close to a non-trivial, contractible 
periodic orbit $x$ of $\rho$. Then 
\begin{equation}
\label{eq:CZ-Delta}
\Delta_\rho(x)-n \leq \CZ(\tx) \leq \Delta_\rho(x)+(n-k)
\end{equation} 
\end{Proposition}

\proof Note that by the continuity of $\Delta$ and
\eqref{eq:mean-index2} we automatically have
$$
\Delta_\rho(x)-n \leq \CZ(\tx) \leq \Delta_\rho(x)+n,
$$
regardless of the nature of the flow of $\rho$. Hence only the second
inequality in \eqref{eq:CZ-Delta} requires a proof.

Furthermore, by arguing as in the proof of Proposition
\ref{prop:index}, it is not hard to reduce the proposition to the
following linear algebra result. Namely, consider a
finite--dimensional symplectic vector space $V$ split as a symplectic
direct sum
$$
V=(L\oplus L^*)\oplus E,
$$ 
where $E$ and $(L\oplus L^*)$ are symplectic spaces, and $L$ and $L^*$ are
Lagrangian in $L\oplus L^*$; cf.\ \eqref{eq:split} and
\eqref{eq:split2}. Set $\dim V=2n$ and $\dim L=k$. Consider a path
$G\colon [0,\,1]\to\Sp(V)$ of the form $G=A\oplus \Gamma$, where
$\Gamma$ is a path in $\Sp(E)$ beginning at $I$ and $A$ is the
block--diagonal path
$$
A=
\begin{bmatrix}
I & tI\\
0 & I
\end{bmatrix},
$$
in $\Sp(L\oplus L^*)$.

\begin{Lemma}
\label{lemma:CZ-Delta}
Let $\tG\colon [0,\,1]\to \Sp(V)$ be a small non-degenerate
perturbation of $G$, also beginning at $I$. Then
\begin{equation}
\label{eq:CZ-Delta2}
\Delta(G)-n\leq \CZ(\tG)\leq \Delta(G)+(n-k)
\end{equation} 
\end{Lemma}

\begin{proof}[Proof of the lemma]
Again, by \eqref{eq:mean-index}, we have 
$$
\Delta(G)-n\leq \CZ(\tG)\leq \Delta(G)+n,
$$
for any path $G$. Hence, only the second
inequality in \eqref{eq:CZ-Delta2} requires a proof. 

Next observe that, once the end-point $\Gamma(1)$ is fixed, the path 
$\Gamma$ is
immaterial for the assertion of the lemma. In other words, if the lemma
holds for one path with a given end-point, it also holds for every
path with the same end-point. This follows from the facts that a
homotopy of $G$ can be traced by a homotopy of $\tG$ (both with fixed
end-points) and that $\CZ$ and $\Delta$ are invariant under such
homotopy and change in the same way when a loop is attached to a path.

As the first step of the proof, let us assume that all eigenvalues of
$\Gamma(1)$ are equal to one. Then $\Gamma(1)$ is in the image of the
exponential mapping $\exp$ for $\Sp(E)$. Indeed, $\Gamma(1)$ is
conjugate to a symplectic linear map which can be chosen to be
arbitrarily close to $I$; see, e.g., \cite[Lemma
5.5]{Gi:conley}. Since $\exp$ is onto a neighborhood of the identity
and commutes with conjugation, $\Gamma(1)$ is in the image of
$\exp$. Furthermore, since $0$ is a regular point of $\exp$ and the
set of regular points is open, we can write $\Gamma(1)=\exp(Q)$, where
$Q$ is a regular point of $\exp$ and all eigenvalues of $Q$ are equal
to zero.

Here we identify the Lie algebra of the symplectic group with the
space of quadratic Hamiltonians. As is customary in symplectic
geometry, the eigenvalues of $Q$ are, by definition, the eigenvalues of
the linear Hamiltonian vector field $X_Q$ generated by $Q$. Also note
that if we identified $\Sp(E)$ with $\Sp(2(n-k))$ and used the matrix
exponential map, we would write $X_Q=JQ$ and $\Gamma(1)=\exp(JQ)$.

Furthermore, we have $A(1)=\exp(\rho)$ in $\Sp(L\oplus L^*)$, where
$\rho$ is a positive definite form on $L^*$ and zero on $L$. Arguing
as above, it is not hard to show that $\rho$ is a regular point of
$\exp$ for $\Sp(L\oplus L^*)$ and that, moreover, $\rho+Q$ is a
regular point of the exponential mapping for $\Sp(V)$. Now we have
$\tG(1)=\exp(K)$ in $\Sp(V)$, where the quadratic form $K$ is close to
$\rho+Q$. In particular, $K$ is also positive definite on $L^*$ and
all eigenvalues of $K$ are close to those of $\rho+Q$, i.e., close to
zero. As has been pointed out above, we can set $\tG(t)=\exp (tK)$ and
$\Gamma(t)=\exp (tQ)$.  As a consequence, with our conventions,
$\CZ(\tG)=-\sgn (K)/2\leq n-k$, where $\sgn (K)$ stands for the
signature of $K$ (i.e., the number of positive eigenvalues minus the
number of negative eigenvalues); see \cite[Section 2.4]{Sa}.  In
addition, $\Delta(G)=0$, and we obtain the second inequality of
\eqref{eq:CZ-Delta2} in this case. To summarize, we have proved
\eqref{eq:CZ-Delta2} when all eigenvalues of $\Gamma(1)$ are equal to
one.

To treat the general case, consider the symplectic direct sum
decomposition $E=E_0\oplus E_1$, where $E_0$ is spanned by the
generalized eigenvectors of $\Gamma(1)$ with eigenvalue one and $E_1$
is the symplectic orthogonal complement of $E_0$ in $E$. Clearly,
$\Gamma(1)$ preserves this decomposition and, after altering if
necessary the path $\Gamma$, we may assume that so do all maps
$\Gamma(t)$.  When $\tG(1)$ is sufficiently close to $G(1)$, we have
the decomposition $V=V_0\oplus V_1$ preserved by $\tG(1)$, where $V_0$
is close to $(L\oplus L^*)\oplus E_0$ and $V_1$ is close to
$E_1$. Applying a time-dependent, close to the identity conjugation to
$\tG(t)$, we reduce the problem to the case where $V_0=(L\oplus
L^*)\oplus E_0$ and $V_1=E_1$. Consider now the paths $G$ and
$\tG$. Both paths begin and end in $\Sp(V_0)\times \Sp(V_1)$, the
first path is contained in this subgroup, and the path $\tG$ is close
to $G$. In particular, $\tG$ is in a tubular neighborhood of the
subgroup. Projecting $\tG$ to $\Sp(V_0)\times \Sp(V_1)$, we can
further reduce the question to the case where $\tG$ is a path in
$\Sp(V_0)\times \Sp(V_1)$, just as $G$ is. Denote by $G=(G_0,G_1)$ and
$\tG=(\tG_0,\tG_1)$ the corresponding decompositions of the paths. The
$E_0$-component of $G_0(1)$ is the map $\Gamma(1)|_{E_0}$ with all
eigenvalues equal to one, and hence \eqref{eq:CZ-Delta2} has already
been proved for $G_0$:
$$
\Delta(G_0)-\dim V_0/2\leq \CZ(\tG_0)\leq \Delta(G_0)+(\dim V_0/2 -k).
$$
On the other hand, the path $\tG_1$ is a small perturbation of the
path $\Gamma|_{E_1}$. Thus, we have
$$
\Delta(G_1)-\dim V_1/2\leq \CZ(\tG_1)\leq \Delta(G_1)+\dim V_1/2.
$$
Recall that $\Delta(G)=\Delta(G_0)+\Delta(G_1)$ and
$\CZ(\tG)=\CZ(\tG_0)+\CZ(\tG_1)$ and that $\dim V_0+\dim V_1=\dim
V=2n$. Thus, adding up these inequalities, we
obtain~\eqref{eq:CZ-Delta2}, which completes the proof of the lemma
and hence the proof of the proposition.
\end{proof}

\subsection{Action selector for ``pinned'' Hamiltonians, following
  E. Kerman}
\label{sec:pinned}
Our goal in this section is to describe a construction of an action
selector for ``pinned'' Hamiltonians, which was introduced in \cite{Ke09,KS}.
Although the class of Hamiltonians and manifolds we work with is
somewhat different from those in \cite{Ke09,KS}, the action selector
is essentially the same as the one considered there. As far as the proofs are
concerned, we adopt here the line of reasoning from \cite{Gi07coiso}
rather than following the Hofer-geometric approach from
\cite{Ke09}. Since the arguments are quite standard, for the sake of
brevity, we just outline the proofs.

Let $M^{2n-k}$ be a closed submanifold, not-necessarily coisotropic,
of a symplectic manifold $W^{2n}$. As before, we
require $W$ to be symplectically aspherical and either closed or a
geometrically bounded and wide.  We assume that $M$ is displaceable
and fix a displaceable open set $U$ containing $M$. Denote by $\HH$
the collection of non-negative, autonomous Hamiltonians $H\colon
W\to\R$ supported in $U$, constant on a small tubular neighborhood of
$M$ and attaining the absolute maximum $C:=\max H$, depending on $H$, on this
neighborhood. Let us require furthermore that $C >e(U)$, where $e(U)$
is the displacement energy of $U$.

It is easy to see that $\HF_n^{(C-\delta,C+\delta)}(H)=\Z_2$ once
$H\in \HH$ and $\delta>0$ is sufficiently small. In fact,
$\HF_*^{(C-\delta,C+\delta)}(H)=\H_{*+n-k}(M;\Z_2)$. Furthermore,
when $a>C$ is large enough (namely, if $a>C+e(U)$), the inclusion map
$$
i_a \colon \Z_2\cong\HF_n^{(C-\delta,C+\delta)}(H)\to \HF_n^{(C-\delta,a)}(H)
$$ 
is zero. The proof of this fact is, for example, contained in the
proof of \cite[Proposition 4.1]{Gi07coiso}; see also \cite{Ke09} for
the case of closed manifolds. This is the main point of the argument
where we need to assume that $W$ is wide (see \cite{Gu07}), unless $W$
is closed. For $H\in \HH$, set
$$
\cf(H)=\inf\{a>C\mid i_a=0\}.
$$
(Strictly speaking, here we have to require $a>C+\delta$ and then also
take infimum over all sufficiently small $\delta>0$.) This is a
version of the action selector for ``pinned'' Hamiltonians, introduced
in \cite{Ke09}.

Alternatively and more explicitly, the action selector $\cf$ can be
defined as follows. Let $\tH$ be a $C^2$-small, non-degenerate
perturbation of $H$, also supported in $U$ (or, to be more precise, in
$S^1\times U$) and such that $\tH\geq H$. Let us also assume that
$\tH$ is autonomous on a small neighborhood of $M$ and that $\max
\tH=C=\max H$ is attained at $p\in M$. (In what follows, we will have
$p$ fixed and independent of $\tH$.)  Then $p$, viewed as an element
of the Floer complex $\CF_*^{(C-\delta,\infty)}(\tH)$, is exact and
there exists a chain in $\CF_{n+1}^{(C-\delta,\infty)}(\tH)$ mapped to
$p$ by the Floer differential; see the proof of \cite[Proposition
4.1]{Gi07coiso}. Let us consider all such chains and, within every
chain, pick an orbit with the largest action and then among the
resulting orbits we choose an orbit $\tx$ with the least action. In
other words, to obtain $\tx$, we first maximize the action within every
chain and then minimize the result among all chains which are
primitives of $p$. Clearly, the orbit $\tx$ is in general not unique,
but the action $\AC_{\tH}(\tx)$ is defined unambiguously. Furthermore,
$\tx$ is connected to $p$ by a Floer downward trajectory and with a
little more effort one can show that in fact $\p \tx=p$. (Again, we
refer the reader to \cite{Gi07coiso} and, in particular, to the proof
of Proposition 4.1 therein for the proofs of these facts; note also
that $\tx$ is the orbit denoted by $\gamma$ in \cite[Proposition
4.1]{Gi07coiso}.) Let us now set $\cf(\tH)=\AC_{\tH}(\tx)$.  Then
$\cf(H)$ is the infimum or the limit (in the obvious sense) of
$\cf(\tH)$ over all such perturbations $\tH$ of $H$. (It is clear that
$\cf(H)$ is less than or equal to the limit; the fact that $\cf(H)$ is
greater than or equal to the limit is a consequence of the definition
of the Floer homology for degenerate Hamiltonians such as $H$.)

It follows from this description that there exists an orbit $x$ of
$H$, referred to in what follows as a \emph{special one-periodic orbit} of
$H$, obtained as a limit point of the orbits $\tx$ in the space of
loops as $\tH\to H$, such that
\begin{equation}
\label{eq:x}
C<\AC_H(x)=\cf(H)< C+e(U)\text{ and } 1\leq\Delta(x)\leq 2n+1.
\end{equation}
Here the first inequalities can be proved by the continuity of the
action (with a little extra argument showing that the inequalities are
strict) and the second ones follow from the continuity of the mean index
and \eqref{eq:mean-index2}. Note that, in general, the special orbit $x$ is
not unique.

We refer the reader to \cite{Ke09} for a detailed investigation of the
properties of the action selector $\cf$. One of these  is particularly
important for our argument.

\begin{Proposition}[\cite{Ke09}]
\label{prop:cont}
The action selector $\cf$ is Lipschitz, with Lipschitz constant equal
to one, on $\HH$ equipped with the $\sup$-norm.
\end{Proposition}

As an immediate consequence of the proposition, the selector $\cf$ extends
from $\HH$ to the $C^0$-closure of $\HH$ in the space of continuous
functions supported in $U$ and this extension is again
Lipschitz with Lipschitz constant equal to one.  For the sake of
completeness, we touch upon a proof of the proposition.

\begin{proof}[Outline of the proof] 
  Let $H$ and $K$ be two Hamiltonians in $\HH$.  Consider the
  perturbations $\tH$ and $\tK$ as above. Clearly, it suffices to show
  that
\begin{equation}
\label{eq:Lips}
|\cf(\tH)-\cf(\tK)|\leq \| \tH-\tK\|_{\Ho},
\end{equation}
where
$$
\|F\|_{\Ho}:=\int_0^1 (\max_W F_t-\min_W F_t)\,dt
$$ 
stands for the Hofer norm of $F$.

Furthermore, denote by $\tx$ again a least action primitive of $p$ in
$\CF_*^{(C-\delta,\infty)}(\tH)$ described above. In particular,
$\cf(\tH)= \AC_{\tH}(\tx)$.  It is not hard to see that under the
linear homotopy from $H$ to $K$, the orbit $\tx$ is mapped to a
primitive $\ty=\sum \ty_i$ of $p$ in the complex
$\CF_*^{(C-\delta,\infty)}(\tK)$, but not necessarily to a least
action primitive. In any case, $\cf(\tK)\leq \AC_{\tK}(\ty):=\max
\AC_{\tK}(\ty_i)$. On the other hand, a standard calculation yields
that
$$
\AC_{\tK}(\ty)-\AC_{\tH}(\tx)\leq \| \tH-\tK\|_{\Ho}.
$$
Hence, we also have $\cf(\tK)-\cf(\tH)|\leq \| \tH-\tK\|_{\Ho}$. A
similar argument, but using the homotopy from $\tK$ to $\tH$, shows
that $\cf(\tH)-\cf(\tK)\leq \| \tH-\tK\|_{\Ho}$, and \eqref{eq:Lips}
follows.
\end{proof}

\begin{Remark}
  It is worth pointing out that the main advantage of using the action
  selector for pinned Hamiltonians in the proof of the main theorem over
  the ordinary action selector is that the former enables us to
  determine the location of the special orbit $x$ via Lemma
  \ref{lemma:location} without additional requirements on $M$ such as
  that $M$ has restricted contact type. This results in sharper index
  and energy bounds that we would have otherwise, cf.\ \cite{Gi07coiso}.
\end{Remark}

\subsection{Proof of Theorem \ref{thm:main}}
\label{subsec:proof}
Throughout the proof, 
as in Section \ref{subsec:prelim}, a neighborhood of $M$ in $W$ is
identified with a neighborhood of $M$ in $M\times \R^k$ equipped with
the symplectic form $\omega = \wM + \sum_{j=1}^k d(p_j \alpha_j)$.
Using this identification, we denote by $U_R$ or just $U$, with $R>0$
sufficiently small, the neighborhood of $M$ in $W$ corresponding to
$M\times B^k_R$. (Thus, $U_R=\{\rho < R^2 /2\}$.) Also set
$|p|:=\sqrt{2\rho}$.

The proof of the theorem relies on a method, by now quite standard, developed in
\cite{Vi90}. The first, albeit technical, step is to specify the class
of ``test'' Hamiltonians.

\subsubsection{The Hamiltonians}
  Fix two real constants $r>0$ and $\eps>0$ with $\eps< r< R$ and a
  constant $C>e(U)$. Let $H\colon [0,\,R]\to \R$ be a smooth, non-negative,
  (non-strictly) decreasing function such that
\begin{itemize}
\item on $[0,\,\eps]$ the function $H$ is a positive constant $C$,
\item on $[\eps,\,2\eps]$ the function $H$ is concave (i.e., $H''\leq 0$),
\item on $[2\eps,\,r-\eps]$ the function $H$ is linear decreasing from
$C-\eps$ to $\eps$,
\item on $[r-\eps,\,r]$ the function $H$ is convex (i.e., $H''\geq 0$),
\item on $[r,\,R]$ the function $H$ is identically zero.
\end{itemize}
Abusing notation, we also denote by $H$ the function equal to $H(|p|)$
on $U$ and equal to zero outside $U$. Let us fix the value of the
parameter $r$, which is not essential for what follows. The parameters
$C$ and $\eps$ will vary and we consider the family of functions
$H=H_{C,\eps}$ parametrized by $C$ and $\eps$ and depending smoothly
on these parameters.

Clearly, $H\in \HH$ for any choice of $\eps$ and $C$.  As $\eps\to 0$,
the functions $H_{C,\eps}$ converge uniformly to the continuous
functions $H_{C,0}$ equal to $C$ on $M$, zero outside $U_r$, and
depending linearly of $|p|$ on $U_r$. It is clear that the limit
functions $H_{C,0}$ are continuous in $C$. Thus, by Proposition
\ref{prop:cont}, $\cf(H_{C,\eps})$ is a continuous function of $C$
and $\eps$ including the limit value $\eps=0$. Moreover, the function
$C\mapsto \cf(H_{C,0})$ is Lipschitz with Lipschitz constant equal to one.

Denote by $X$ the Hamiltonian vector field of the function $|p|$ on
$U\setminus M$.  By Proposition \ref{prop:flow}, the integral curves of
$X$ project to the geodesics of the leaf-wise metric $\rho$ on $M$,
parametrized by arc length. The Hamiltonian vector field of $H$ is
$$
X_H=H'X,
$$
where $H'$ stands for the derivative of $H$ with respect to
$|p|$. Note that even though $X$ is defined only on $U\setminus M$, the vector
field $X_H$ is defined everywhere, for $H$ is constant near $M$ and
outside $U_r$. Thus, non-trivial one-periodic orbits of $X_H$ lie on
the levels $|p|=\const$ with $H'(|p|)$ in the length spectrum $\SC$ of
the metric $\rho$. (Recall that, by definition, $\SC$ is formed by the
lengths of non-trivial closed leaf-wise geodesics of $\rho$. Here, we
may restrict our attention only to the geodesics contractible in $W$.)
Furthermore, observe that the ``coordinates'' $p_i$ are constant along
the orbits of the flow of $X_H$. In other words, every trajectory
starting in $U$ lies on a coisotropic submanifold $M\times p\subset
U$.  This is a particular case of conservation of momentum.

Let $x$ be a non-trivial one-periodic orbit of $H$. A direct
calculation relying on Proposition \ref{prop:normalform} shows that
\begin{eqnarray*}
\AC_H(x)&=&H(x)+\AC(x)\\
&=&H(x)+\AC(\pi(x))-|p(x)| l(\pi(x)),
\end{eqnarray*}
where $l$ and $\AC$ stand for the length of the curve and,
respectively, the negative symplectic area bounded by the curve.

Assume that the slope of $H$ (on the interval $[2\eps,\,r-\eps]$)
is outside $\SC$. (This is a generic condition.) Then the orbit $x$
lies on the level where $|p(x)|$ is either in the range
$[\eps,\,2\eps]$ or in the range $[r-\eps,r]$.  Let now $x$ be a
special one-periodic orbit from Section \ref{sec:pinned} such that, in
particular, \eqref{eq:x} holds.  The key to the proof is the following
lemma, which specifies the location of $x$ for, at least, some sequence of
the Hamiltonians $H$.

\begin{Lemma}
\label{lemma:location}
There exists a sequence $C_j\to\infty$  such that the slopes of all
functions $H_{C_j,\eps}$, with $\eps>0$ sufficiently small, are outside $\SC$ 
and $|p(x)|\in [\eps,\,2\eps]$.
\end{Lemma}

In the Lagrangian case this observation can be traced back to the
original work of Viterbo, \cite{Vi90}. Here we follow the treatment
from \cite{Ke09} with several modifications resulting from our
somewhat different conventions and more importantly from the fact that
$M$ is now coisotropic.

\begin{proof}[Proof of Lemma \ref{lemma:location}] 
  The slope of the function $H_{C,0}$ is $C/r$. This slope is in $\SC$
  if and only if $C\in r\SC$ in the obvious notation. The set $\SC$
  (and hence $r\SC$) is closed, and the slope of $H_{C,\eps}$ is close
  to the slope of $H_{C,0}$ when $\eps>0$ is small. As a consequence,
  the slope of $H_{C,\eps}$ is outside $\SC$ whenever $C\not\in r\SC$
  and $\eps>0$ is small.

  Pick $C\not\in r\SC$ and a positive sequence $\eps_i\to 0$.  Without
  loss of generality, we may require all $\eps_i$ to be sufficiently
  close to zero to ensure that the slope of $H_i$ is not in $\SC$.
  Let $x_i$ be a special orbit of $H_i:=H_{C,\eps_i}$.  Since the
  norms of the differentials $dH_i$ are bounded from above, the norms
  of the derivatives $\dot{x}_i$ are point-wise bounded.  By the
  Arzela--Ascoli theorem, we may assume, after passing if necessary to
  a subsequence, that the orbits $x_i$ converge to a curve $y$ lying
  on a level $|p|=\const$ including possibly the submanifold $M$. It
  is clear that $y$ is smooth and projects to a closed, leaf-wise
  geodesic on $M$. Furthermore,
  $$
  \AC_{H_i}(x_i)=\cf(H_i)\to H_{C,0}(y)+\AC(y)=\cf(H_{C,0}),
  $$
  by the continuity of the action functional and of the action selector
  $\cf$.

  If $|p(x_i)|$ is in the range $[r-\eps_i,\,r]$ for all $i$, the
  orbit $y$ is on the level $|p|=r$ and $H_{C,0}(y)=0$. Thus, we then have
\begin{equation}
\label{eq:Sigma}
\cf(H_{C,0})=\AC(y)\in \Sigma,
\end{equation}
where $\Sigma$ is the action spectrum or, to be more precise, the
symplectic area spectrum of the level $|p|=r$, i.e., the collection of
symplectic areas bounded by contractible closed characteristics on
this level. 

Arguing by contradiction, assume now that the lemma fails, i.e., for
every sufficiently large $C$, say $C>a$, which is not in $r\SC$, there
exists such a sequence $\eps_i$ with $|p(x_i)|$ in the range
$[r-\eps_i,\,r]$. Consider the function $f(C):=\cf(H_{C,0})$ on the
interval $[a,\,\infty)$.  By \eqref{eq:Sigma}, $f$ sends the set
$[a,\,\infty)\setminus r\SC$ to $\Sigma$. Recall that $r\SC$ is
not only closed, but also has zero measure; see \cite[Lemma
6.6]{Gi07coiso}. Furthermore, by Proposition \ref{prop:cont}, $f$ is a
Lipschitz function and, as is well known (see, e.g., \cite{HZ94}),
$\Sigma$ has measure zero. To summarize, $f$ is a Lipschitz function
sending a full measure set to a zero measure set. Such a function is
necessarily constant. This is impossible, for $f(C)\geq C$
by~\eqref{eq:x}.
\end{proof}

Let us fix one of the constants $C=C_j$ from Lemma
\ref{lemma:location} and let $H_i=H_{C_j,\eps_i}$. Denote by $x_i$, or
just $x$, its one-periodic orbit such as in the lemma.  (For the proof
of the theorem we do not need the entire double sequence, but only one
family of Hamiltonians $H_{C_j,\eps_i}$ parametrized by $\eps_i$.)
Clearly, $\gamma_i=\pi(x_i)$ is a leaf-wise geodesic on $M$. Since the
slopes of Hamiltonians $H_i$ are bounded from above (by, say,
$2C_j/r$), it is easy to prove using the Arzela--Ascoli theorem that
the geodesics $\gamma_i$ converge as $i\to\infty$ after if necessary
passing to a subsequence. Denote the limit geodesic (traversed in the
opposite direction) by $\eta$. Our goal is to show that $\eta$ has the
required properties \eqref{eq:bound-index} and \eqref{eq:bound-area}.
The fact that, by Lemma \ref{lemma:location}, $|p(x_i)|\in
[\eps_i,\,2\eps_i]$ (i.e., $x_i$ lies in the region where $H_i$ is
concave) will be essential for proving this.

\subsubsection{Index bounds}
Consider a perturbation $\tH$ of $H=H_i$ as in Section
\ref{sec:pinned}. This Hamiltonian has a one-periodic orbit $\tx$, a
perturbation of $x=x_i$, with index $n+1$. After reparametrizing $x$
and reversing its orientation, we can view $x$ as a periodic orbit
$x^-$ of $\rho$. Likewise, $\tx$ can be viewed as a periodic orbit
$\tx^-$ of a non-degenerate perturbation $K$ of $\rho$. Denote by
$\gamma^-=\pi(x^-)$ the geodesic $\gamma=\gamma_i$ with reversed
orientation.

By Proposition \ref{prop:index}, we have
$$
\mu(\gamma^-)=-\Delta_\rho(x^-),
$$
and thus, by Proposition \eqref{prop:CZ-Delta},
$$
-\mu(\gamma^-)-n\leq\CZ(\tx^-)\leq -\mu(\gamma^-)+ (n-k).
$$
It is not hard to show that $\CZ(\tx^-)=-\CZ(\tx)=-(n+1)$ using the
fact that $x$ is in the region where $H$ is concave (i.e., $|p(x)|\in
[\eps_i,\,2\eps_i]$) by Lemma \ref{lemma:location}.  As a consequence,
$$
n+1 \leq \mu(\gamma^-)+n \text{ and } \mu(\gamma^-)-n+k\leq n+1.
$$
Hence,
$$
1\leq \mu(\gamma^-)\leq 2n+1-k.
$$
Passing to the limit and using the continuity of the mean index, we
conclude that the same holds for $\eta$, the limit of the curves
$\gamma^-$. This proves \eqref{eq:bound-index}.

\begin{Remark}
  Note that if we have used here just the second inequality of
  \eqref{eq:x} rather than Proposition \eqref{prop:CZ-Delta}, we would
  have the weaker bound $1\leq \mu(\gamma^-)\leq
  2n+1$. 
\end{Remark}

\subsubsection{Action bounds}
By the first inequality in \eqref{eq:x}, we have
\begin{equation}
\label{eq:aux}
C<\AC_H(x)=H(x)+\AC(\gamma)-|p(x)|l(\gamma)<C+e(U).
\end{equation}
Here, by the definition of $H$ and Lemma \ref{lemma:location},
$|p(x)|\in [\eps_i,\,2\eps_i]$ and $H(x)\in [C,C-\eps_i]$. Note that
the sequence $l(\gamma)$ with $\gamma=\gamma_i$ is bounded as
$i\to\infty$ due the fact that the slope of $H_i$ is bounded. Thus,
passing to the limit (for a subsequence if necessary), we have $0\leq
-\AC(\eta)\leq e(U)$. Here, the negative sign comes from the fact that
$\eta$ is the limit of $\gamma^-$, i.e., the geodesics $\gamma$ with
reversed orientation. Taking $r>0$ sufficiently small, we obtain
$$
0\leq \Area(\eta)\leq e(M)+\delta
$$
for any given $\delta>0$, where $\Area(\eta)=-\AC(\eta)$ is the symplectic area
bounded by $\eta$. To finish the proof, we need to ensure that the
first inequality is strict: $\Area(\eta)>0$. This is an immediate
consequence of the non-trivial fact that, by \cite[Theorem
6.1]{Gi07coiso}, $\AC_H(x)-C\geq \eps$ for some $\eps>0$
independent of $i$. For then, by the first inequality in
\eqref{eq:aux}, $\Area(\gamma^-)>\eps/2$ when $i$ is large
enough. This concludes the proof of \eqref{eq:bound-area}, and thus
the proof of the theorem.



\begin{thebibliography}{GGW}

\bibitem[AF1]{AF08}
P. Albers, U. Frauenfelder,
Leaf-wise intersections and Rabinowitz Floer homology, Preprint 2008,
arXiv:0810.3845.


\bibitem[AF2]{AF08g}
P. Albers, U. Frauenfelder,
Infinitely many leaf-wise intersection points on cotangent bundles,
Preprint 2008, arXiv:0812.4426.


\bibitem[ALP]{AL}
  M. Audin, F. Lalonde, L. Polterovich, Lagrangian submanifolds, in
  \emph{Holomorphic Curves in Symplectic Geometry}, Eds: M. Audin,
  J. Lafontaine,  Progress in Mathematics, \textbf{117}, Birkh\"auser
  Verlag, Basel, 1994, pp.\ 271--321.


\bibitem[Ba]{Ba}
A. Banyaga,
On fixed points of symplectic maps, \emph{Invent.\ Math.}, \textbf{56}
(1980), 215--229.


\bibitem[Bo1]{Bo96}
P. Bolle, 
Une condition de contact pour les sous-vari\'et\'es co\"isotropes
d'une vari\'et\'e symplectique, \emph{C.\ R.\ Acad.\ Sci.\ Paris,
S\'erie I}, \textbf{322} (1996), 83--86.


\bibitem[Bo2]{Bo98}
P. Bolle, 
A contact condition for p-dimensional submanifolds of a symplectic
manifold ($2\leq p\leq n$), \emph{Math.\ Z.}, \textbf{227} (1998),
211--230.

\bibitem[CGK]{CGK}
K. Cieliebak, V. Ginzburg, E. Kerman,
Symplectic homology and periodic orbits near symplectic submanifolds,
\emph{Comment.\ Math.\ Helv.}, \textbf{79} (2004), 554--581.


\bibitem[Dr]{Dr}
D. Dragnev, 
Symplectic rigidity, symplectic fixed points and global perturbations
of Hamiltonian systems, \emph{Comm.\ Pure Appl.\ Math.} \textbf{61}
(2008), 346--370.

\bibitem[Du]{Du}
J.J. Duistermaat, 
On the Morse index in variational calculus, \emph{Adv.\ Math.}, \textbf{21}
(1984), 207-253.


\bibitem[EH]{EH} 
I. Ekeland, H. Hofer, 
Two symplectic fixed-point theorems with applications to Hamiltonian
dynamics, \emph{J.\ Math.\ Pures Appl.} \textbf{68} (1989), 467--489.

\bibitem[FHW]{FHW}
A. Floer, H. Hofer, K. Wysocki,
Applications of symplectic homology, I. \emph{Math.\ Z.}, \textbf{217} (1994),
577--606.

\bibitem[Gi1]{Gi99:orbits} V.L. Ginzburg, Hamiltonian dynamical
  systems without periodic orbits, in \emph{Northern California
    Symplectic Geometry Seminar}, Eds: Ya. Eliashberg, D. Fuchs, T. Ratiu,
  A. Weinstein, Amer.\ Math.\ Soc.\ Transl.\ Ser.\ 2, vol.\
  196, Amer.\ Math.\ Soc., Providence, RI, 1999, pp.\ 35--48.

\bibitem[Gi2]{Gi:alan} 
V.L. Ginzburg, The Weinstein conjecture and the theorems of nearby and 
almost existence,  in \emph{The Breadth of Symplectic and Poisson Geometry.
Festschrift in Honor of Alan Weinstein}; J.E. Marsden and T.S. Ratiu
(Eds.), Birkh\"auser, 2005, pp.\ 139--172.

\bibitem[Gi3]{Gi07coiso}
V.L. Ginzburg,
Coisotropic intersections, \emph{Duke Math.\ J.}, \textbf{140} (2007),
111--163.

\bibitem[Gi4]{Gi:conley}
V.L. Ginzburg, The Conley conjecture, Preprint 2006, math.SG/0610956; to appear in
\emph{Ann.\ of Math.}


\bibitem[GG1]{GG03:seifert}
V.L. Ginzburg, B.Z. G\"urel,
A $C^2$-smooth counterexample to the Hamiltonian Seifert 
conjecture in $\R^4$, \emph{Ann.\ of Math.}, \textbf{158} (2003), 953--976.

\bibitem[GG2]{GG04}
V.L. Ginzburg, B.Z. G\"urel,
Relative Hofer--Zehnder capacity and periodic orbits in twisted
cotangent bundles, \emph{Duke Math.\ J.}, \textbf{123} (2004), 1--47.


\bibitem[GG3]{GG09wm}
V.L. Ginzburg, B.Z. G\"urel,
Periodic orbits of twisted geodesic flows and the Weinstein-Moser theorem,  
\emph{Comment.\ Math.\ Helv.},  \textbf{84}  (2009), 865-907. 

\bibitem[G\"u1]{Gu07}
B.Z. G\"urel, 
Totally non-coisotropic displacement and its applications to
Hamiltonian dynamics, \emph{Comm.\ Contemp.\ Math.} \textbf{10}
(2008),  1103--1128.


\bibitem[G\"u2]{Gu09imrn}
B.Z. G\"urel, Leafwise coisotropic intersections, 
\emph{Int.\ Math.\ Res.\ Not.\ IMRN},  2009; doi: 10.1093/imrn/rnp164.

\bibitem[G\"u3]{Gu:ex}
B.Z. G\"urel, 
Fragility of leafwise coisotropic intersections, in preparation.


\bibitem[Ho]{Ho90a}
H. Hofer, 
On the topological properties of symplectic maps, \emph{Proc.\ Roy.\
Soc.\ Edinburgh Sect.\ A}, \textbf{115} (1990),  25--38.


\bibitem[HZ]{HZ94}
H. Hofer, E. Zehnder,
\emph{Symplectic Invariants and Hamiltonian Dynamics}, Birk\"auser,
1994.

\bibitem[Ka]{Ka}
J. Kang, Existence of leafwise intersection points in the unrestricted case,
Preprint 2009, arXiv:0910.2369.


\bibitem[Ke1]{Ke07coiso} 
E. Kerman,
Displacement energy of coisotropic submanifolds and Hofer's geometry,
\emph{J.\ Mod.\ Dyn.} \textbf{2} (2008), 471--497.

\bibitem[Ke2]{Ke09}
E. Kerman,
Action selectors and Maslov class rigidity, 
\emph{Int.\ Math.\ Res.\ Not.\ IMRN},  2009; doi:10.1093/imrn/rnp093.

\bibitem[KS]{KS}
E. Kerman, N.I. Sirikci,
Maslov class rigidity for Lagrangian submanifolds via Hofer's geometry. 
Preprint 2008, arXiv:0808.1422. To appear in  \emph{Comment.\ Math.\ Helv.}.

\bibitem[Lo]{Lo}
Y. Long,  \emph{Index Theory for Symplectic Paths with Applications}, 
Progress in Mathematics, 207. Birkh\"auser Verlag, Basel, 2002. 

\bibitem[Mo]{Mo}
J. Moser, 
A fixed point theorem in symplectic geometry, \emph{Acta Math.},
\textbf{141} (1978), 17--34.

\bibitem[Oh]{Oh}
Y.-G. Oh,
Geometry of coisotropic submanifolds in symplectic and K\"ahler manifolds,
Preprint 2003, math.SG/0310482. 

\bibitem[Po1]{Po91Tr} L. Polterovich, The Maslov class of the Lagrange
  surfaces and Gromov's pseudo-holomorphic curves,
  \emph{Trans.\ Amer.\ Math.\ Soc.},  \textbf{325} (1991), 241--248.

\bibitem[Po2]{Po91MZ} L. Polterovich, Monotone Lagrange submanifolds
  of linear spaces and the Maslov class in cotangent bundles,
  \emph{Math.\ Z.}, \textbf{207} (1991), 217--222.

\bibitem[Sa]{Sa}
D.A. Salamon, Lectures on Floer homology, in
\emph{Symplectic Geometry and Topology},  Eds: Y. Eliashberg and
L. Traynor, IAS/Park City Mathematics series, \textbf{7}, 1999,
pp.\ 143--230.

\bibitem[SZ]{SZ}
D. Salamon, E. Zehnder, 
Morse theory for periodic solutions of Hamiltonian systems and the
Maslov index, \emph{Comm.\ Pure Appl.\ Math.},\textbf{45} (1992),
1303--1360.

\bibitem[Si]{Si}
J.-C. Sikorav,
Some properties of holomorphic curves in almost complex manifolds,
in
  \emph{Holomorphic Curves in Symplectic Geometry}, Eds: M. Audin,
  J. Lafontaine, Progress in Mathematics, \textbf{117}, Birkh\"auser
  Verlag, Basel, 1994, pp.\ 165--189.

\bibitem[Su]{Su} 
D. Sullivan, A foliation of geodesics is characterized
by having no ``tangent homologies'', \emph{J. Pure Appl.\ Algebra},
\textbf{13} (1978), 101--104.

\bibitem[To]{To}
B. Tonnelier,
in preparation.


\bibitem[Us]{U09coiso}
M. Usher, 
Boundary depth in Floer theory and its applications to Hamiltonian
dynamics and coisotropic submanifolds, Preprint 2009, arXiv:0903.0903; to appear in
\emph{Israel J. Math.}

\bibitem[Vi]{Vi90}
C. Viterbo, 
A new obstruction to embedding Lagrangian tori,
\emph{Invent.\ Math.}, \textbf{100} (1990),  301--320. 

\bibitem[We]{We}
J. Weber,
Perturbed closed geodesics are periodic orbits: Index and transversality,
\emph{Math.\ Z.}, \textbf{241} (2002), 45--81.


\bibitem[Zi1]{Zilt08}
F. Ziltener,
Coisotropic Submanifolds, leaf wise fixed points, and symplectic
embeddings, Preprint 2008, arXiv:0811.3715v2. To appear in the \emph{J. Sympl.\
Geom.}

\bibitem[Zi2]{ZiltPrep} F. Ziltener, A Maslov map for coisotropic
  submanifolds, leaf-wise fixed points and presymplectic
  non-embeddings, Preprint 2009, arXiv:0911.1460

\end{thebibliography}
\end{document}